\documentclass[a4paper,12pt]{amsart}
\usepackage{amssymb, amsmath}
\usepackage[curve]{xypic}

\providecommand{\Ker}{\textnormal{Ker}}
\providecommand{\IIm}{\textnormal{Im}}

\providecommand{\Ker}{\textnormal{Ker}}

\providecommand{\Int}{\textnormal{int}}

\providecommand{\Hol}{\textnormal{Hol}}
\providecommand{\Tor}{\textnormal{Tor}}

\providecommand{\fl}{\textnormal{fl}}
\providecommand{\CS}{\textnormal{CS}}
\providecommand{\Int}{\textnormal{int}}

\begin{document}

\title{Relative Deligne cohomology and Cheeger-Simons characters}
\author{Fabio Ferrari Ruffino}
\address{ICMC - Universidade de S\~ao Paulo, Avenida Trabalhador s\~ao-carlense 400, 13566-590 - S\~ao Carlos - SP, Brasil}
\email{ferrariruffino@gmail.com}
\thanks{The author was supported by FAPESP (Funda\c{c}\~ao de Amparo \`a Pesquisa do Estado de S\~ao Paulo).}

\begin{abstract}
In \cite{BT} the authors discuss two possible definitions of the relative Cheeger-Simons characters, the second one fitting into a long exact sequence. Here we relate that picture to the one of the relative Deligne cohomology groups, defined via the mapping cone: we show that there are three meaningful relative groups, and we analyze the corresponding definition of relative Cheeger-Simons characters in each case. We then extract another definition, corresponding to the one fitting into a long exact sequence in \cite{BT}. Finally we show how the explicit formulas for the holonomy, the transgression maps and the integration can be extended to the relative case.
\end{abstract}

\maketitle

\newtheorem{Theorem}{Theorem}[section]
\newtheorem{Lemma}[Theorem]{Lemma}
\newtheorem{Corollary}[Theorem]{Corollary}
\newtheorem{Rmk}[Theorem]{Remark}
\newtheorem{Def}{Definition}[section]
\newtheorem{ThmDef}[Theorem]{Theorem - Definition}

\section{Introduction}

A Deligne cohomology class \cite{Brylinski} of degree $p$ on a smooth manifold $X$ is a hypercohomology class of degree $p$ of the complex of sheaves:
\begin{equation}\label{DeligneComplex}
	\underline{U}(1) \overset{\tilde{d}}\longrightarrow \Omega^{1}_{\mathbb{R}} \overset{d}\longrightarrow \cdots \overset{d}\longrightarrow \Omega^{p}_{\mathbb{R}},
\end{equation}
where $\underline{U}(1)$ is the sheaf of smooth $U(1)$-valued functions, $\Omega^{k}_{\mathbb{R}}$ is the sheaf of smooth real differential forms of degree $k$ and $\tilde{d} = \frac{1}{2\pi i} d \circ \log$. The curvature is a globally defined closed differential form of degree $p+1$. Given such a class it is possible to define its holonomy as a $U(1)$-valued function defined on smooth singular $p$-cycles \cite{GT2}, and, when the curvature vanishes, the holonomy only depends on the corresponding homology class. Moreover, if $M$ is a smooth manifold of dimension $m$ and $\mathcal{M} = C^{\infty}(M, X)$, there are transgression maps:
	\[\psi_{M}: H^{p}(X, \underline{U}(1) \overset{\tilde{d}}\longrightarrow \cdots \overset{d}\longrightarrow \Omega^{p}_{\mathbb{R}}) \longrightarrow H^{p-m}(\mathcal{M}, \underline{U}(1) \overset{\tilde{d}}\longrightarrow \cdots \overset{d}\longrightarrow \Omega^{p-m}_{\mathbb{R}}).
\]
When $m = p$ the map $\psi_{M}$ provides a function $\mathcal{M} \rightarrow U(1)$ which, on a point $f \in \mathcal{M}$, coincides with the holonomy on the cycle $f_{*}[M]$. Finally, if $E \rightarrow X$ is a smooth oriented fiber bundle of rank $k$, there is an integration map in Deligne cohomology \cite{GT}:
	\[\int_{E/X}: H^{p}(E, \underline{U}(1) \overset{\tilde{d}}\longrightarrow \cdots \overset{d}\longrightarrow \Omega^{p}_{\mathbb{R}}) \longrightarrow H^{p-k}(X, \underline{U}(1) \overset{\tilde{d}}\longrightarrow \cdots \overset{d}\longrightarrow \Omega^{p-k}_{\mathbb{R}}).
\]
The Deligne cohomology group of degree $p$ is canonically isomorphic to the group of Cheeger-Simons differential characters of the same degree. An element of the latter is a couple $(\chi, \omega)$, where $\chi$ is a $U(1)$-valued group morphism defined on the smooth $p$-cycles of $X$ (the holonomy), and $\omega$ is an integral $(p+1)$-form on $X$ (the curvature) such that, on a $p$-boundary $\partial D$, one has:
	\[\chi(\partial D) = \int_{D} \omega \mod \mathbb{Z}.
\]

Relative versions of Deligne cohomology and Cheeger-Simons characters have been studied by several authors \cite{BT, Zucchini, HL, Shahbazi}. In \cite{Zucchini} the author describes a relative version of Deligne cohomology via concrete computations in $\rm\check{C}$ech cohomology, he also provides a definition of relative Cheeger-Simons characters and shows some important physical applications. In \cite{HL} the authors consider differential characters relative to the boundary of a manifold in order to formulate the Lefschetz-Pontrjagin Duality Theorem in this context. In the present paper we refer in particular to \cite{BT}, where the authors show how to define in two possible ways the relative version of the Cheeger-Simons characters, the second one fitting into a long exact sequence. Moreover, we consider the possible relative versions of the Deligne cohomology, and we describe explicitly the relations between these two pictures \cite{Shahbazi, FR}: it turns out that there are three meaningful relative Deligne cohomology groups, and we analyze the corresponding definitions of relative Cheeger-Simons character, all of which are directly related to the first definition in \cite{BT}. We show that these groups actually fit into long exact sequences, but mixing the topological cohomology groups and their differential extensions. The definition leading to the long exact sequence in \cite{BT} can also be extracted from the Deligne cohomology groups, even if in a less natural way. Moreover, there are explicit formulas for the holonomy, the transgression maps and the integration in the absolute case \cite{GT2, GT}, and we show how to extend them to the relative case in the three possible formulations. In particular, the relative formula for the holonomy is used to prove the isomorphisms between the relative Deligne cohomology groups and the corresponding groups of relative Cheeger-Simons characters.

The paper is organized as follows. In section \ref{RelDeligne} we recall the definition of relative Deligne cohomology and we show the corresponding long exact sequence in cohomology in each case. In section \ref{RelCS} we show the possible definitions of the relative Cheeger-Simons characters. In section \ref{Holonomy} we provide the formula for the the holonomy and the transgression maps in the relative case, and we show the formula for the integration map when the fiber has a boundary.

\section{Relative Deligne cohomology}\label{RelDeligne}

We first recall some basic facts about trivializations in Deligne cohomology, then we discuss the relative Deligne cohomology groups.

\subsection{Deligne cohomology and trivializations}

Let us consider the complex \eqref{DeligneComplex} on a manifold $X$, which we call $S^{\bullet}_{X, p}$. We call $C^{n}(S^{\bullet}_{X, p})$, $Z^{n}(S^{\bullet}_{X, p})$ and $B^{n}(S^{\bullet}_{X, p})$ the groups of cochains, cocycles and coboundaries of degree $n$. We have recalled in the introduction the main features of the hypercohomology group of degree $p$, i.e.\ $H^{p}(S^{\bullet}_{X, p})$. Moreover, $H^{n}(S^{\bullet}_{X, p}) \simeq H^{n}(X, \mathbb{R}/\mathbb{Z})$ for $n < p$ and $H^{n}(S^{\bullet}_{X, p}) \simeq H^{n+1}(X, \mathbb{Z})$ for $n > p$ \cite{Brylinski}. There is a natural projection:
	\[\pi_{\underline{U}(1)}: C^{n}(S^{\bullet}_{X, p}) \rightarrow C^{n}(X, \underline{U}(1)).
\]
For $n = p$, the first Chern class of $x \in H^{p}(S^{\bullet}_{X, p})$ is the image of $[\pi_{\underline{U}(1)}\tilde{x}]$ under the isomorphism $H^{p}(X, \underline{U}(1)) \simeq H^{p+1}(X, \mathbb{Z})$, where $\tilde{x} \in Z^{p}(S^{\bullet}_{X, p})$ is any cocycle representing $x$.
\begin{Def}\label{Trivializations} Given a cocycle $\alpha \in Z^{p}(S^{\bullet}_{X, p})$, we call it:
\begin{itemize}
	\item \emph{geometrically trivial} if $\alpha \in B^{p}(S^{\bullet}_{X, p})$; a \emph{geometric trivialization} of $\alpha$ is a class of cochains $[\beta] \in C^{p-1}(S^{\bullet}_{X, p})/B^{p-1}(S^{\bullet}_{X, p})$ such that $\delta^{p-1}\beta = \alpha$;
	\item \emph{topologically trivial} if $[\pi_{\underline{U}(1)}\alpha] = 0$; a \emph{topological trivialization} of $\alpha$ is a class of cochains $[\beta'] \in C^{p-1}(X, \underline{U}(1))/B^{p-1}(X, \underline{U}(1))$ such that $\delta^{p-1}\beta' = \pi_{\underline{U}(1)}\alpha$.
\end{itemize}
Moreover:
\begin{itemize}
	\item if $\alpha$ is topologically trivial, a \emph{strong topological trivialization} of $\alpha$ is a class of cochains $[\beta''] \in C^{p-1}(S^{\bullet}_{X, p})/B^{p-1}(S^{\bullet}_{X, p})$ such that $\alpha \cdot (\delta^{p-1}\beta'')^{-1} = (1, 0, \ldots, 0, \rho)$, for $\rho$ a $p$-form on $X$.
\end{itemize}
\end{Def}
The geometric trivializations of $\alpha$ are a coset of the group of the geometric trivializations of $(1, 0, \ldots, 0)$, which is $H^{p-1}(S^{\bullet}_{X, p}) \simeq H^{p-1}(X, \mathbb{R}/\mathbb{Z}) \simeq H^{p-1}(X, U(1))$. Similarly, the topological trivializations of $\alpha$ are a coset of $H^{p-1}(X, \underline{U}(1))$. For $\alpha$ geometrically trivial, there is a natural map from the geometric trivializations to the topological ones (even on the representatives): if $\alpha = (1, 0, \ldots, 0)$ the image is the subset of the topological trivializations admitting a locally constant representative. For $p > 1$ the latter are the torsion ones, and the kernel of the map is $H^{p-1}(X, \mathbb{R})/H^{p-1}(X, \mathbb{Z})$.\footnote{For $p = 1$, this corresponds to the fact that a topological trivialization of a line bundle is a global section, while a geometric one is a parallel global section. In this case the map is injective, in particular the kernel is not $H^{p-1}(X, \mathbb{R})/H^{p-1}(X, \mathbb{Z})$ for $p = 1$. In fact, if we consider the factorization:
	\[H^{p-1}(X, \mathbb{R}/\mathbb{Z}) \overset{a}\longrightarrow H^{p-1}(X, \underline{U}(1)) \overset{b} \longrightarrow H^{p}(X, \mathbb{Z}),
\]
one has $\Ker(b \circ a) \simeq H^{p-1}(X, \mathbb{R})/H^{p-1}(X, \mathbb{Z})$, and, for $p > 1$, $\Ker(a) = \Ker(b \circ a)$ since $b$ is an isomorphism. For $p = 1$, instead, $b$ is not injective, because its kernel contains the functions admitting a global logarithm. In this case $a$ is injective (and it is the map we are interested to), and $b \circ a = 0$, since $H^{0}(X, \mathbb{R}/\mathbb{Z}) \simeq H^{0}(X, \mathbb{R})/H^{0}(X, \mathbb{Z})$ or, equivalently, $\Tor \, H^{1}(X, \mathbb{Z}) = 0$.}

The group of the strong topological trivializations of $(1, 0, \ldots, 0)$ is naturally isomorphic to $H^{p-1}(S^{\bullet}_{X, p-1})$: in fact, $C^{p-1}(S^{\bullet}_{X, p}) = C^{p-1}(S^{\bullet}_{X, p-1})$ and $B^{p-1}(S^{\bullet}_{X, p}) = B^{p-1}(S^{\bullet}_{X, p-1})$, and the condition for $\beta''$ representing a strong topological trivialization is the cocycle condition in $C^{p-1}(S^{\bullet}_{X, p-1})$. There is a natural surjective map from the strong topological trivializations to the topological trivializations, which, for $\alpha = (1, 0, \ldots, 0)$, coincides with natural morphism $H^{p-1}(S^{\bullet}_{X, p-1}) \rightarrow H^{p-1}(X, \underline{U}(1))$. When $\alpha$ is geometrically trivial, the geometric trivializations are a subset of the strong topological trivializations (in this case the form $\rho$ in def.\ \ref{Trivializations} is closed and integral for a topological trivialization, while it is $0$ for a geometric trivialization).

\subsection{Relative Deligne cohomology}

Given a smooth map $f: Y \rightarrow X$ of manifolds without boundary, we consider the complex \eqref{DeligneComplex} in degree $p$ on $X$ and $q$ on $Y$, which we call:
\begin{equation}\label{DeligneComplexXY}
	S^{\bullet}_{X,p} := \underline{U}(1)_{X} \overset{\tilde{d}}\longrightarrow \cdots \overset{d}\longrightarrow \Omega^{p}_{X,\mathbb{R}} \qquad S^{\bullet}_{Y,q} := \underline{U}(1)_{Y} \overset{\tilde{d}}\longrightarrow \cdots \overset{d}\longrightarrow \Omega^{q}_{Y,\mathbb{R}}.
\end{equation}
We consider the push-forward $f_{*}S^{\bullet}_{Y,q}$, recalling that, for $\mathcal{F}$ a sheaf on $Y$, the sheaf $f_{*}\mathcal{F}$ on $X$ is defined as $(f_{*}\mathcal{F})(U) := \mathcal{F}(f^{-1}U)$ for any $U \subset X$ open. For $q \leq p$ there is a natural map of complexes of sheaves on $X$:
\begin{equation}\label{MapVarphi}
	\varphi_{f,p,q}: S^{\bullet}_{X,p} \rightarrow f_{*}S^{\bullet}_{Y,q},
\end{equation}
defined in the following way: in degree $i \leq q$ it pulls back via $f$ from $U$ to $f^{-1}U$ the function or differential form in the domain, in degree $i > q$ it is the zero-map.\footnote{It is necessary that $q \leq p$, otherwise \eqref{MapVarphi} would not be a map of complexes, since it would not commute with the coboundary of degree $p$.} One can construct the mapping cone of the morphism \eqref{MapVarphi}, i.e.\ the complex of sheaves on $X$ \cite{GM}:\footnote{The matrix defining the coboundary is supposed to multiply from the left a column vector.}
\begin{equation}\label{Cone}
	C(\varphi_{f,p,q})^{\bullet} := S^{\bullet}_{X,p} \oplus f_{*}S^{\bullet-1}_{Y,q} \qquad d^{\bullet}_{C(\varphi_{X,Y,p,q})} := \begin{bmatrix} d^{\bullet}_{S_{X,p}} & 0 \\ \varphi_{f,p,q}^{\bullet} & -f_{*}d^{\bullet-1}_{S_{Y,q}} \end{bmatrix}.
\end{equation}
\begin{Def}\label{Relative} The \emph{relative Deligne cohomology groups of type $(p,q)$} of the map $f: Y \rightarrow X$ are the hypercohomology groups of \eqref{Cone}. We denote them as $H^{n}(S_{X,p}, S_{Y,q}, f)$, and we denote the groups of cocycles and coboundaries as $Z^{n}(S_{X,p}, S_{Y,q}, f)$ and $B^{n}(S_{X,p}, S_{Y,q}, f)$.
\end{Def}
For $p = q$ the long exact sequence in cohomology corresponds, up to canonical isomorphism, to the following sequence (we denote by $H^{\bullet}(X, Y, f; G)$ the relative cohomology groups with coefficients in $G$ associated to the map $f$):
\begin{equation}\label{LongDeligne1}
\begin{split}
	\cdots & \longrightarrow H^{p-1}(X, Y, f; \mathbb{R}/\mathbb{Z}) \longrightarrow H^{p-1}(X; \mathbb{R}/\mathbb{Z}) \longrightarrow H^{p-1}(Y; \mathbb{R}/\mathbb{Z}) \\
	& \longrightarrow H^{p}(S_{X,p}, S_{Y,p}, f) \longrightarrow H^{p}(S_{X,p}) \longrightarrow H^{p}(S_{Y,p}) \\
	& \longrightarrow H^{p+2}(X, Y, f; \mathbb{Z}) \longrightarrow H^{p+2}(X; \mathbb{Z}) \longrightarrow H^{p+2}(Y; \mathbb{Z}) \longrightarrow \cdots.
\end{split}
\end{equation}
A class in $H^{p}(S_{X,p}, S_{Y,p}, f)$ is represented by a cocycle in $\tilde{\alpha} \in Z^{p}(S_{X,p})$ together with a representative of a geometric trivialization of $f^{*}\tilde{\alpha}$ on $Y$. The curvature is therefore a closed $(p+1)$-form on $X$ which vanishes when pulled back to $Y$. A class, whose projection in $H^{p}(S_{X,p})$ is trivial, corresponds by exactness to a geometrical trivialization of $(1, 0, \ldots, 0)$ on $Y$, up to the pull-back of a geometric trivialization on $X$. Therefore, a coboundary is an element of $\tilde{\alpha} \in B^{p}(S_{X,p})$ together with a representative of a geometric trivialization of $f^{*}\tilde{\alpha}$, which, up to $(p-2)$-coboundaries, is the pull-back of a representative of a geometric trivialization of $\tilde{\alpha}$.\footnote{In particular, a coboundary in $B^{p}(S_{X,p})$ can be extended to a non-trivial relative class, if it is endowed with a trivialization on $Y$ which is not the pull-back of any one on $X$.} The subgroup of flat classes in $H^{p}(S_{X,p}, S_{Y,p}, f)$ corresponds to the relative cohomology group $H^{p}(X, Y, f; \mathbb{R}/\mathbb{Z})$. \\

For $q = p-1$ the long exact sequence in cohomology corresponds, up to canonical isomorphism, to the following sequence:
\begin{equation}\label{LongDeligne2}
\begin{split}
	\cdots & \longrightarrow H^{p-1}(X, Y, f; \mathbb{R}/\mathbb{Z}) \longrightarrow H^{p-1}(X; \mathbb{R}/\mathbb{Z}) \longrightarrow H^{p-1}(S_{Y,p-1}) \\
	& \longrightarrow H^{p}(S_{X,p}, S_{Y,p-1}, f) \longrightarrow H^{p}(S_{X,p}) \longrightarrow H^{p+1}(Y; \mathbb{Z}) \\
	& \longrightarrow H^{p+2}(X, Y, f; \mathbb{Z}) \longrightarrow H^{p+2}(X; \mathbb{Z}) \longrightarrow H^{p+2}(Y; \mathbb{Z}) \longrightarrow \cdots.
\end{split}
\end{equation}
A class in $H^{p}(S_{X,p}, S_{Y,p-1})$ is represented by a cocycle in $\tilde{\alpha} \in Z^{p}(S_{X,p})$ together with a representative of a strong topological trivialization of $f^{*}\tilde{\alpha}$ on $Y$. The curvature is therefore a closed $(p+1)$-form on $X$ which is exact when pulled back to $Y$. A class, whose projection in $H^{p}(S_{X,p})$ is trivial, corresponds by exactness to a strong topological trivialization of $(1, 0, \ldots, 0)$ on $Y$, up to the pull-back of a \emph{geometric} trivialization on $X$. Therefore, the coboundaries coincide with the ones in $B^{p}(S_{X,p}, S_{Y,p}, f)$.\footnote{In particular, a coboundary in $B^{p}(S_{X,p})$ can be extended to a non-trivial relative class, if it is endowed with a strong topological trivialization on $Y$ that is not geometrical, or if it is endowed with a geometrical trivialization which is not the pull-back of any one on $X$.} In particular, there is a natural embedding:
\begin{equation}\label{MapPP1}
	\varphi_{f,p,p-1}: H^{p}(S_{X,p}, S_{Y,p}, f) \hookrightarrow H^{p}(S_{X,p}, S_{Y,p-1}, f),
\end{equation}
induced from the embedding $Z^{p}(S_{X,p}, S_{Y,p}, f) \hookrightarrow Z^{p}(S_{X,p}, S_{Y,p-1}, f)$. More precisely, there is a morphism of exact sequences from \eqref{LongDeligne1} to \eqref{LongDeligne2}, since the map $H^{p-1}(Y; \mathbb{R}/\mathbb{Z}) \rightarrow H^{p-1}(S_{Y,p-1})$ corresponds to the embedding of the flat classes, and the surjective map $H^{p}(S_{Y,p}) \rightarrow H^{p+1}(Y; \mathbb{Z})$ is the first Chern class. Such a map is induced from the natural map of complexes $(S_{X,p} \rightarrow S_{Y,p}) \rightarrow (S_{X,p} \rightarrow S_{Y,p-1})$. \\

\paragraph{\emph{Remark}:}\label{Intermediate} We have seen that an element of $\alpha \in H^{p}(S_{X,p}, S_{Y,p-1}, f)$ is represented by a cocycle $\tilde{\alpha} \in Z^{p}(S_{X,p})$ together with a representative of a strong topological trivialization of $f^{*}\tilde{\alpha}$. The fact of $f^{*}\tilde{\alpha}$ being geometrically trivial is not enough to conclude that $\alpha$ belongs to the image of the embedding \eqref{MapPP1}, since the topological trivialization is not necessarily a geometric one. In particular, when $f^{*}\tilde{\alpha}$ is geometrically trivial, the form $\rho$ appearing in the definition \ref{Trivializations} is integral, but not necessarily $0$. Therefore, the subgroup of $H^{p}(S_{X,p}, S_{Y,p-1}, f)$ containing the classes that are geometrically trivial on $Y$ properly contains the image of \eqref{MapPP1}. $\square$ \\

For $q \leq p-2$ the long exact sequence in cohomology, starting from two positions before $H^{p}(S_{X,p}, S_{Y,q}, f)$, always corresponds, up to canonical isomorphism, to the following sequence:
\begin{equation}\label{LongDeligne3}
\begin{split}
	\phantom{\cdots} & \phantom{XXXXX} \cdots \phantom{XXX} \longrightarrow H^{p-1}(X; \mathbb{R}/\mathbb{Z}) \longrightarrow H^{p}(Y; \mathbb{Z}) \\
	& \longrightarrow H^{p}(S_{X,p}, S_{Y,q}, f) \longrightarrow H^{p}(S_{X,p}) \longrightarrow H^{p+1}(Y; \mathbb{Z}) \\
	& \longrightarrow H^{p+2}(X, Y, f; \mathbb{Z}) \longrightarrow H^{p+2}(X; \mathbb{Z}) \longrightarrow H^{p+2}(Y; \mathbb{Z}) \longrightarrow \cdots.
\end{split}
\end{equation}
We start for simplicity from the case $q = 0$. A class in $H^{p}(S_{X,p}, S_{Y,0}, f)$ is represented by a cocycle in $\tilde{\alpha} \in Z^{p}(S_{X,p})$ together with a representative of a topological trivialization of $f^{*}\tilde{\alpha}$ on $Y$. The curvature is therefore a closed $(p+1)$-form on $X$ which is exact when pulled back to $Y$, as for the case $q = p-1$. A class, whose projection in $H^{p}(S_{X,p})$ is trivial, corresponds by exactness to a topological trivialization of $(1, 0, \ldots, 0)$ on $Y$, up to the pull-back of a \emph{geometric} trivialization on $X$, composed with the natural map from geometric to topological trivializations. Therefore, a coboundary is an element of $\tilde{\alpha} \in B^{p}(S_{X,p})$ together with a representative of a topological trivialization of $f^{*}\tilde{\alpha}$, which restricts, up to $(p-2)$-coboundaries, the pull-back of a representative of a geometric trivialization of $\tilde{\alpha}$. When $q > 0$, the description is analogous, but all the representatives of the topological trivializations involved are extended to cochains of degree $p-1$ in $S_{Y,q}$; this has no effect on the cohomology classes involved. There is a natural surjective map:
\begin{equation}\label{MapP10}
	\varphi_{f,p-1,0}: H^{p}(S_{X,p}, S_{Y,p-1}, f) \rightarrow H^{p}(S_{X,p}, S_{Y,0}, f),
\end{equation}
induced from the surjective map $Z^{p}(S_{X,p}, S_{Y,p-1}, f) \rightarrow Z^{p}(S_{X,p}, S_{Y,0}, f)$. More precisely, there is a morphism of exact sequences from \eqref{LongDeligne2} to \eqref{LongDeligne3}, since the map $H^{p-1}(S_{Y,p-1}) \rightarrow H^{p}(Y; \mathbb{Z})$ is the first Chern class. Such a map is induced from the natural map of complexes $(S_{X,p} \rightarrow S_{Y,p-1}) \rightarrow (S_{X,p} \rightarrow S_{Y,0})$.

For what concerns the first part of the long exact sequence, it actually depends on $q$. In particular, from the degree $0$ to $q-1$, it corresponds to the sequence in singular cohomology with $\mathbb{R}/\mathbb{Z}$-coefficients. Even $H^{q}(S_{X,p}, S_{Y,q}, f) \simeq H^{q}(X, Y, f; \mathbb{R}/\mathbb{Z})$, as one can prove from the definition or using the five lemma in the exact sequence, hence we arrive at $H^{q}(X, \mathbb{R}/\mathbb{Z})$. Then there are the five terms $H^{q}(X, \mathbb{R}/\mathbb{Z}) \rightarrow H^{q}(S_{Y,q}) \rightarrow H^{q+1}(S_{X,p}, S_{Y,q}, f) \rightarrow H^{q+1}(X, \mathbb{R}/\mathbb{Z}) \rightarrow H^{q+2}(Y, \mathbb{Z})$. Thus, $H^{q+1}(S_{X,p}, S_{Y,q}, f)$ is the subgroup of $H^{q+1}(S_{X,q+1}, S_{Y,q}, f)$ made by classes which are flat on $X$. Then, for $q < n < p-1$, there are blocks of the form $H^{n}(X, \mathbb{R}/\mathbb{Z}) \rightarrow H^{n+1}(Y, \mathbb{Z}) \rightarrow H^{n+1}(S_{X,p}, S_{Y,q}, f) \rightarrow H^{n+1}(X, \mathbb{R}/\mathbb{Z}) \rightarrow H^{n+2}(Y, \mathbb{Z})$. Thus, $H^{n+1}(S_{X,p}, S_{Y,q}, f)$ is the subgroup of $H^{n+1}(S_{X,n+1}, S_{Y,0}, f)$ made by classes which are flat on $X$. \\

We can also consider another definition of the relative Deligne cohomology groups, which seems more artificial within this picture, but has the advantage of fitting into a long exact sequence all made by Deligne cohomology groups. We consider the following maps:
\begin{itemize}
	\item $H^{p}(S_{X,p}, S_{Y,p-1}, f) \rightarrow H^{p}(S_{Y,p})$, defined composing the map $H^{p}(S_{X,p}, S_{Y,p-1},$ $f) \rightarrow H^{p}(S_{X,p})$ appearing in the sequence \eqref{LongDeligne2} with the restriction map $H^{p}(S_{X,p}) \rightarrow H^{p}(S_{Y,p})$ appearing in the sequence \eqref{LongDeligne1};
	\item $H^{p-1}(S_{X,p-1}) \rightarrow H^{p}(S_{X,p}, S_{Y,p-1}, f)$, defined composing the restriction map $H^{p-1}(S_{X,p-1}) \rightarrow H^{p-1}(S_{Y,p-1})$ appearing in the sequence \eqref{LongDeligne1} with the Bockstein map $H^{p-1}(S_{Y,p-1}) \rightarrow H^{p}(S_{X,p}, S_{Y,p-1}, f)$ appearing in the sequence \eqref{LongDeligne2}.
\end{itemize}
We define:
\begin{equation}\label{HBar}
	\overline{H}^{p}(S_{X,p}, S_{Y,p}, f) := \frac{\Ker(H^{p}(S_{X,p}, S_{Y,p-1}, f) \rightarrow H^{p}(S_{Y,p}))}{\IIm(H^{p-1}(S_{X,p-1}) \rightarrow H^{p}(S_{X,p}, S_{Y,p-1}, f))}.
\end{equation}
Clearly for $Y = \emptyset$ we get $\overline{H}^{p}(S_{X,p}) \simeq H^{p}(S_{X,p})$.
\begin{Theorem} There is a long exact sequence:
\begin{equation}\label{LongDeligne4}
	\cdots \longrightarrow H^{p-1}(S_{Y,p-1}) \overset{\beta}\longrightarrow \overline{H}^{p}(S_{X,p}, S_{Y,p}, f) \longrightarrow H^{p}(S_{X,p}) \longrightarrow H^{p}(S_{Y,p}) \longrightarrow \cdots.
\end{equation}
\end{Theorem}
\paragraph{Proof:} The Bockstein map $\beta$ is induced from the one of \eqref{LongDeligne2}, that we call $\beta'$: the image of $\beta'$ is contained in the numerator of \eqref{HBar} because of the exactness of \eqref{LongDeligne2}. By definition of the denominator of \eqref{HBar} the kernel of $\beta$ is the image of the restriction map $H^{p-1}(S_{X,p-1}) \rightarrow H^{p-1}(S_{Y,p-1})$, thus \eqref{LongDeligne4} is exact in $H^{p-1}(S_{Y,p-1})$. The image of $\beta'$ contains the denominator of \eqref{HBar}, therefore the exactness in $\overline{H}^{p}(S_{X,p}, S_{Y,p}, f)$ follows from the one of \eqref{LongDeligne2}. Finally, in order to prove the exactness in $H^{p}(S_{X,p})$, we consider the following commutative diagram:
	\[\xymatrix{
	H^{p}(S_{X,p}, S_{Y,p}, f) \ar[r]^(.6){\eta} \ar[d]_{\psi} & H^{p}(S_{X,p}) \\
	\overline{H}^{p}(S_{X,p}, S_{Y,p}, f) \ar[ur]_(.55){\nu}
}\]
where $\eta$ is the map appearing in \eqref{LongDeligne1}, $\nu$ the one appearing in \eqref{LongDeligne4} and $\psi$ is the composition of the embedding $H^{p}(S_{X,p}, S_{Y,p}, f) \rightarrow H^{p}(S_{X,p}, S_{Y,p-1}, f)$, whose image is contained in the numerator of \eqref{HBar} because of the exactness of \eqref{LongDeligne1}, with the projection to the quotient in \eqref{HBar}. We show that $\IIm \, \eta = \IIm \, \nu$, so that the exactness of \eqref{LongDeligne4} in $H^{p}(S_{X,p})$ follows from the one of \eqref{LongDeligne1}. Obviously $\IIm \, \eta \subset \IIm \, \nu$. For the converse, the image of the embedding $H^{p}(S_{X,p}, S_{Y,p}, f) \rightarrow H^{p}(S_{X,p}, S_{Y,p-1}, f)$ is the subset of classes which are trivial when pulled-back to $H^{p}(S_{Y,p})$, therefore, applying $\nu$ to the numerator of \eqref{HBar}, we get classes belonging to the kernel of $H^{p}(S_{X,p}) \rightarrow H^{p}(S_{Y,p})$, i.e.\ to the image of $\eta$. $\square$

\section{Cheeger-Simons characters}\label{RelCS}

We fix some notations as in \cite{BT}. For a given map of manifolds $f: Y \rightarrow X$:
\begin{itemize}
	\item we call $C_{\bullet}(f)$ the chain complex $C_{\bullet}(X) \times C_{\bullet-1}(Y)$ with boundary $\partial(C, C') = (\partial C + f_{*}C', -\partial C')$; we call $Z_{\bullet}(f)$ and $B_{\bullet}(f)$ the subgroups of cycles and boundaries;
	\item we call $\Omega^{\bullet}(f)$ the cochain complex $\Omega^{\bullet}(X) \times \Omega^{\bullet-1}(Y)$ with coboundary $d(\omega, \rho) = (d\omega, f^{*}\omega - d\rho)$; we call $\Omega^{\bullet}_{cl}(f)$ the subgroup of relative closed forms, i.e.\ such that $d(\omega, \rho) = 0$;
	\item we call $\Omega^{\bullet}_{0}(f)$ the complex of smooth forms $\omega$ on $X$ that $f^{*}\omega = 0$.
\end{itemize}
A couple $(\omega, \rho) \in \Omega^{\bullet}_{cl}(f)$ represents a relative de-Rham cohomology class via the integration $\int_{(C,C')}(\omega, \rho) := \int_{C} \omega + \int_{C'}\rho$. Because of the natural embedding $\Omega^{\bullet}_{0}(f) \hookrightarrow \Omega^{\bullet}(f)$ defined by $\omega \rightarrow (\omega, 0)$, the same holds for a form $\omega \in \Omega^{\bullet}_{0, cl}(f)$. Let us consider $\omega \in \Omega^{p+1}_{cl}(X)$ representing an integral class on $X$, and such that $f^{*}\omega$ is exact, i.e.\ such that the de-Rahm class of $\omega$ on $X$ is liftable to a class relative to $Y$ via $f$. We call $\Lambda_{\omega} \subset \Omega^{p}(Y)$ the subset containing the forms $\rho$ such that $(\omega, \rho)$ represents an integral relative class (in particular $d\rho = f^{*}\omega$). For $\omega = 0$, $\Lambda_{\omega}$ is the subgroup of closed forms on $Y$ which are integral on the kernel of the push-forward $f_{*}: H_{p}(Y) \rightarrow H_{p}(X)$. In fact, considering $(C, C') \in Z_{p+1}(f)$, one has that $\rho \in \Lambda_{\omega}$ if and only if $\int_{C'} \rho = -\int_{C} \omega + n$ for $n \in \mathbb{Z}$, and, because of the cycle condition, $C'$ is a cycle in $Y$ such that $[f_{*}C'] = 0$. Given a class $\alpha \in H^{p}(S_{X,p}, S_{Y,p-1}, f)$, if $\omega$ is the curvature and $(1, 0, \ldots, 0, \rho)$ is the reparametrization induced on $Y$ by the strong topological trivialization, then $\rho \in \Lambda_{\omega}$.\footnote{In order to prove the last statement, let us consider a relative cycle $(C,C')$. Then the holonomy on $C'$ of the image of $\alpha$ via the map $H^{p}(S_{X,p}, S_{Y,p-1}, f) \rightarrow H^{p}(S_{X,p}) \rightarrow H^{p-1}(S_{Y,p-1})$ can be computed as $\exp 2\pi i \int_{C'}\rho$ or as $\exp 2\pi i \int_{-C} \omega$, therefore $\int_{C} \omega + \int_{C'} \rho \in \mathbb{Z}$.} \\

From now on we suppose that all of the singular chains are smooth. We provide a definition of relative Cheeger-Simons character for each of the three meaningful relative Deligne cohomology groups introduced above.
\begin{Def} For $f: Y \rightarrow X$ a map of manifolds, a \emph{relative Cheeger-Simons character of type I} of degree $p$ is a couple $(\chi, \omega)$ where:
\begin{itemize}
	\item $\chi: Z_{p}(f) \rightarrow \mathbb{R}/\mathbb{Z}$ is a group homomorphism;
	\item $\omega \in \Omega^{p+1}_{0, cl}(f)$ represents an integral relative de-Rham cohomology class;
	\item given $(C, C') \in B_{p}(f)$ such that $(C, C') = \partial(D, D')$, one has:
	\begin{equation}\label{RelativeCS1}
	\chi(C, C') = \int_{D} \omega \mod \mathbb{Z}.
\end{equation}
\end{itemize}
\end{Def}
Formula \eqref{RelativeCS1} is consistent, since the result does not depend on the choice of $(D, D')$. In fact, if we choose another trivialization $(E, E')$, then \eqref{RelativeCS1} gives the same result if and only if $\int_{D - E} \omega \in \mathbb{Z}$. Since $\partial(D - E, D' - E') = 0$ and $(\omega, 0)$ represents an integral relative cohomology class, the thesis follows.

The next definition coincides with \cite[Definition 2.1]{BT}.
\begin{Def} If $(X, Y)$ is a pair of manifolds, a \emph{relative Cheeger-Simons character of type II} of degree $p$ is a triple $(\chi, \omega, \rho)$ where:
\begin{itemize}
	\item $\chi: Z_{p}(f) \rightarrow \mathbb{R}/\mathbb{Z}$ is a group homomorphism;
	\item $(\omega, \rho) \in \Omega^{p+1}_{cl}(f)$ represents an integral relative de-Rham cohomology class;
	\item given $(C, C') \in B_{p}(f)$ such that $(C, C') = \partial(D, D')$, one has:
	\begin{equation}\label{RelativeCS2}
	\chi(C, C') = \int_{D} \omega + \int_{D'} \rho \mod \mathbb{Z}.
\end{equation}
\end{itemize}
\end{Def}
Formula \eqref{RelativeCS2} is consistent because of the same argument used after formula \eqref{RelativeCS1}. We now consider the trivial principal $\mathbb{R}/\mathbb{Z}$-bundle $\mathcal{L} \rightarrow Z_{p-1}(Y)$ whose fiber over $A \in Z_{p-1}(Y)$ is the set of equivalence classes of triples $[(\xi, \eta, t)]$, where $(\xi, \eta)$ is a Cheeger-Simons character of degree $p-1$ on $Y$, $t \in \mathbb{R}/\mathbb{Z}$, and $(\xi, \eta, t) \sim (\xi + \xi', \eta + \eta', t - \xi'(A))$ for every $(\xi', \eta')$. Every character $(\xi, \eta)$ defines a global section $[(\xi, \eta, 0)]$, and, because of the global section $[(0, 0, 0)]$, the bundle is canonically trivial.
\begin{Def} If $(X, Y)$ is a pair of manifolds, a \emph{relative Cheeger-Simons character of type III} of degree $p$ is an equivalence class of couples $[(\chi, \omega)]$ where:
\begin{itemize}
	\item $\omega \in \Omega^{p+1}(X)$ is integral and $\omega\vert_{Y}$ is exact;
	\item $\chi \in \Gamma(\pi_{2}^{*}\mathcal{L})$, for $\pi_{2}: Z_{p}(f) \rightarrow Z_{p-1}(Y)$ the projection on the second term;
	\item there exists $\rho \in \Lambda_{\omega}$ such that, given $(C, C') \in B_{p}(f)$ verifying $(C, C') = \partial(D, D')$, one has:
	\begin{equation}\label{RelativeCS3}
	\chi(C, C')/\pi_{2}^{*}([(0, 0, 0)]_{C'}) = \int_{D} \omega + \int_{D'} \rho \mod \mathbb{Z}.
\end{equation}
\end{itemize}
Two couples $(\chi, \omega)$ and $(\chi', \omega)$ are equivalent if there exists an automorphisms $\varphi: \mathcal{L} \rightarrow \mathcal{L}$ of the form $\varphi[(\xi, \eta, 0)] = [(\xi + \xi', \eta + \eta', 0)]$ such that $\pi_{2}^{*}\varphi \circ \chi = \chi'$.
\end{Def}
Given a couple $(\chi, \omega)$ representing a type III relative character, we get a type II relative character $(\tilde{\chi}, \omega, \rho)$ where $\tilde{\chi}(C, C') := \chi(C, C')/\pi_{2}^{*}([(0, 0, 0)]_{C'})$ and $\rho$ is defined by equation \eqref{RelativeCS3}. Fixing a character $(\xi, \eta)$ on $Y$ we determine a trivialization $[(\xi, \eta, 0)]$ of $\mathcal{L}$, and we obtain another type II relative character $(\tilde{\chi}', \omega, \rho')$ where $\tilde{\chi}'(C, C') := \chi(C, C')/\pi_{2}^{*}([(\xi, \eta, 0)]_{C'}) = \tilde{\chi}(C, C') - \xi(C')$ and, from equation \eqref{RelativeCS3}, $\rho' = \rho + \eta$. Equivalently, if we consider the automorphism $\varphi[(\xi', \eta', 0)] = [(\xi' + \xi, \eta' + \eta, 0)]$ of $\mathcal{L}$, then $\tilde{\chi}'(C, C') := (\pi_{2}^{*}\varphi \circ \chi)(C, C')/\pi_{2}^{*}([(0, 0, 0)]_{C'})$. This corresponds to a natural action of the (absolute) characters of $Y$ on the group of relative characters of type II of $(X,Y)$, defined as follows. For $(\chi, \omega, \rho)$ a type II character on $(X,Y)$ and $(\xi,\eta)$ a character on $Y$:
\begin{equation}\label{ActionTypeII}
	((\xi, \eta) \cdot \chi)(C, C') := \chi(C) - \xi(C') \qquad (\xi, \eta) \cdot \omega := \omega \qquad (\xi, \eta) \cdot \rho := \rho + \eta.
\end{equation}
The group of the relative characters of type III of $(X, Y)$ is therefore the orbit space of the action \eqref{ActionTypeII}. \\

\paragraph{}An explicit isomorphism from $H^{p}(S_{X,p}, S_{Y,p-t}, f)$, for $t = 0, 1, 2$, to the group of relative Cheeger-Simons $p$-characters of type $t+1$ can be obtained via the holonomy of a relative Deligne class, for which we provide the explicit formula in the following. Their is a natural immersion from characters of type I to characters of type II defined by $(\chi, \omega) \rightarrow (\chi, \omega, 0)$, which corresponds to the embedding \eqref{MapPP1}. Moreover, their is a natural surjective map from characters of type II to characters of type III, that corresponds to the surjective map \eqref{MapP10}, defined by \eqref{RelativeCS3} for a fixed character $(\chi, \omega, \rho)$. We can actually be more precise. We call $\CS^{p}_{t}(f)$ the group of relative characters of degree $p$ and type $t$, for $t = I, II, III$. Considering the remark a few lines after formula \eqref{MapPP1}, the embedding $\CS^{p}_{I}(f) \hookrightarrow \CS^{p}_{II}(f)$ factorizes through an intermediate step $\CS^{p}_{I}(f) \hookrightarrow \CS^{p}_{II'}(f) \hookrightarrow \CS^{p}_{II}(f)$, where $\CS^{p}_{II'}(f)$ is the group of type II relative characters $(\chi, \omega, \rho)$ whose pull-back to $Y$ is geometrically trivial, i.e.\ such that $\rho$ is integral \cite{BT}. Moreover, calling $\Omega^{\bullet+1}_{int}(f)$ the group of relative closed forms representing an integral class, we define:
	\[\Omega^{\bullet}_{f}(Y) := \IIm(\pi_{2}: \Omega^{\bullet+1}_{int}(f) \rightarrow \Omega^{\bullet}(Y)).
\]
In other words, $\rho \in \Omega^{\bullet}_{f}(Y)$ if and only if it can be completed to a type II relative character $(\chi, \omega, \rho)$. If $f$ is an embedding then $\Omega^{\bullet}_{f}(Y) = \Omega^{\bullet}(Y)$, because, given $\rho \in \Omega^{p-1}(Y)$, one can extend it to $\tilde{\rho} \in \Omega^{p-1}(X)$ and consider the character $(\chi, d\tilde{\rho}, \rho)$ for $\chi(C, C') := \int_{C} \tilde{\rho}$. In general, because of the same argument, $\Omega^{\bullet}_{f}(Y)$ contains all the forms $\rho$ on $Y$ which are pull-back via $f$ of a form $\tilde{\rho}$ on $X$, but it can be larger than this space.
\begin{Theorem} There are three exact sequences fitting into the following commutative diagram:
	\[\xymatrix{
	0 \ar[r] & \CS^{p}_{I}(f) \ar[r]^{i_{1}} \ar@{_(->}[d]^(.45){\iota_{1}} & \CS_{II}^{p}(f) \ar[r]^{\pi_{1}} \ar@{=}[d] & \Omega^{p-1}_{f}(Y) \ar[r] \ar@{->>}[d]^(.4){p_{1}} & 0\, \\
	0 \ar[r] & \CS^{p}_{II'}(f) \ar[r]^{i_{2}} & \CS_{II}^{p}(f) \ar[r]^(.35){\pi_{2}} & \Omega^{p-1}_{f}(Y)/\Omega^{p-1}_{\Int}(Y) \ar[r] & 0\, \\
	0 \ar[r] & \CS^{p-1}(Y)/\CS^{p-1}_{\fl}(X) \ar[r]^(.65){i_{3}} \ar@{^(->}[u]_(.45){\iota_{2}} & \CS_{II}^{p}(f) \ar[r]^{\pi_{3}} \ar@{=}[u] & \CS^{p}_{III}(f) \ar[r] \ar@{->>}[u]_(.4){p_{2}} & 0,
}\]
where $\CS^{p}_{\fl}(X)$ is the group of flat characters on $X$. Of course we get an equivalent commutative diagram replacing the groups of Cheeger-Simons characters with the corresponding relative Deligne cohomology groups.
\end{Theorem}
\paragraph{Proof:} For the first line, we define $i_{1}(\chi, \omega) := (\chi, \omega, 0)$ and $\pi_{1}(\chi, \omega, \rho) := \rho$. The surjectivity of $\pi_{1}$ follows from the definition of $\Omega^{p-1}_{f}(Y)$. For the second line, we define $i_{2}(\chi, \omega, \rho) := (\chi, \omega, \rho)$ and $\pi_{2}(\chi, \omega, \rho) := [\rho]$, and the exactness is obvious. For the third line, we define $i_{3}[(\xi, \eta)] := (\chi, 0, \eta)$ for $\chi(C, C') := -\xi(C')$. It is the Bockstein map of the sequence \eqref{LongDeligne2} from degree $p-1$ to degree $p$ up to its kernel, and this proves the injectivity. We define $\pi_{3}(\chi, \omega, \rho) := [(\chi, \omega)]$ identifying $\chi$ with $\chi/\pi_{2}^{*}[(0, 0, 0)]$ in equation \eqref{RelativeCS3}. The exactness of the sequence follows from the action \eqref{ActionTypeII}, since such an action is exactly the addition of $i_{3}[(\xi, \eta)]$. The map $\iota_{1}$ and $\iota_{2}$ are defined respectively as $i_{1}$ and $i_{3}$, simply restricting the codomain (the fact that the image of $i_{3}$ is contained in $\CS^{p}_{II'}(f)$ can be shown directly or via the exactness of \eqref{LongDeligne2}). Finally, $p_{1}$ is the obvious projection and $p_{2}[(\xi, \omega)] := [\rho]$, for $\rho$ defined by equation \eqref{RelativeCS3}. $\square$ \\

One could get other exact sequences considering the subgroups of characters which are flat on $X$ or only on $Y$ in each case. We skip the details since they are not particularly enlightening, except for the exact sequences in \cite[Theorem 2.4]{BT}.

\subsection{Forth definition.} In \cite[Chapter 4]{BT} the authors define the relative Cheeger-Simons characters also from the relative cochain complex defined by Hopkins and Singer, and they get groups that fit into a long exact sequence. We recall that definition, considering that the group $\CS^{p}_{II'}(f)$ corresponds to the group $\hat{H}^{k}_{0}(\rho)$ defined in \cite{BT}. There is a natural embedding:
\begin{equation}\label{phif}
	\phi_{f}: \Omega^{p-1}_{int}(M)/\Omega^{p-1}_{0,int}(f) \rightarrow \CS^{p}_{II'}(f)
\end{equation}
defined by $\phi_{f}([\rho])(C, C') := \int_{C} \rho$. Then:
\begin{Def} The groups of \emph{relative Cheeger-Simons characters of type IV} is the quotient $\CS^{p}_{II'}(f)/\IIm\,\phi_{f}$.
\end{Def}
We show in the following (theorem \ref{ThmXip3}) the natural isomorphism $\CS^{p}_{IV}(f) \simeq \overline{H}^{p}(S_{X,p}, S_{Y,p}, f)$, the latter being defined in \eqref{HBar}. It follows that these characters fit into a long exact sequence analogous to \eqref{LongDeligne4}, which coincides with the one considered in \cite{BT}.

\section{Transgression maps}\label{Holonomy}

Let $X$ be a compact manifold without boundary and $M$ a compact oriented manifold of dimension $m$ even with boundary. We define:
	\[\mathcal{M}_{X} := C^{\infty}(M, X) \qquad \mathcal{M}_{\partial} := C^{\infty}(\partial M, X).
\]
Given a cover $\mathfrak{U}$ of $X$, in \cite[Definition 2.1]{GT2} the authors define a cover $\mathfrak{V}$ of $\mathcal{M}$ naturally associated to $\mathfrak{U}$, and a transgression map
\begin{equation}\label{Transgression}
	\psi_{\mathfrak{U}}: C^{k}(\mathfrak{U}, S_{X,p}) \rightarrow C^{k-m}(\mathfrak{V}, S_{\mathcal{M}_{X},p-m}).
\end{equation}
Let $r: \mathcal{M} \rightarrow \mathcal{M}_{\partial}$ be the restriction map, and $(\partial \psi)_{\mathfrak{U}}$ the map \eqref{Transgression} obtained considering $\partial M$ instead of $M$. If $D_{\mathfrak{U}}$ and $D_{\mathfrak{V}}$ are the differentials of the associated double complexes, then \eqref{Transgression} satisfies \cite[Theorem 2.1]{GT2}
\begin{equation}\label{TransgressionBoundary}
	\psi_{\mathfrak{U}} \circ D_{\mathfrak{U}} = (-1)^{m}D_{\mathfrak{V}} \circ \psi_{\mathfrak{U}} + r^{*} \circ (\partial \psi)_{\mathfrak{U}}.
\end{equation}
Let us now consider a map of manifolds $f: Y \rightarrow X$. We define:
	\[\mathcal{M}_{Y} := C^{\infty}(\partial M, Y) \qquad \mathcal{M}_{f} = \{(g, g') \in \mathcal{M}_{X} \times \mathcal{M}_{Y}: g\vert_{\partial M} = f \circ g'\}.
\]
There is a commutative diagram:
\begin{equation}\label{DiagramM}
	\xymatrix{
	& \mathcal{M}_{X} \ar[dr]^{r} & \\
	\mathcal{M}_{f} \ar[ur]^{\pi_{1}} \ar[dr]_{\pi_{2}} & & \mathcal{M}_{\partial} \\
	& \mathcal{M}_{Y} \ar[ur]_{l_{f}}
}
\end{equation}
where $l_{f}(g') = f \circ g'$.

We call $C^{k}(\mathfrak{U}, S_{X,p}, S_{Y,p-1}, f)$ the group of cochains of the complex \eqref{Cone} for $q = p-1$ associated to $\mathfrak{U}$ in degree $k$. An element $\alpha \in C^{k}(\mathfrak{U}, S_{X,p}, S_{Y,p-1}, f)$ is a couple $(\alpha_{1}, f_{*}\alpha_{2})$ with $\alpha_{1} \in C^{k}(\mathfrak{U}, S_{X,p})$ and $\alpha_{2} \in C^{k-1}(f^{-1}\mathfrak{U}, S_{Y,p-1})$.
\begin{Def} Let $f: Y \rightarrow X$ be a map of manifolds without boundary and $M$ a compact manifold of dimension $m$ even with boundary. For $\mathfrak{V}$ the cover induced on $\mathcal{M}_{X}$ by $\mathfrak{U}$ and $\mathfrak{V}'$ the one induced on $\mathcal{M}_{Y}$ by $f^{-1}\mathfrak{U}$, we put $\mathfrak{W} = (\pi_{1}, \pi_{2})^{-1}(\mathfrak{V} \times \mathfrak{V}')$ on $\mathcal{M}_{f}$. We define the \emph{relative transgression map}
\begin{equation}\label{RelativeTransgression}
	\psi_{\mathfrak{U},f}: C^{k}(\mathfrak{U}, S_{X,p}, S_{Y,p-1}, f) \rightarrow C^{k-m}(\mathfrak{W}, S_{\mathcal{M}_{f},p-m})
\end{equation}
as:
	\[\psi_{\mathfrak{U}, f}(\alpha_{1}, f_{*}\alpha_{2}) := \pi_{1}^{*}\psi_{\mathfrak{U}}(\alpha_{1}) - \pi_{2}^{*}\psi_{f^{-1}\mathfrak{U}}(\alpha_{2}).
\]
\end{Def}
\begin{Lemma}\label{RelativeTransgressionBoundaryThm} If $D_{\mathfrak{U}}$ and $D_{\mathfrak{W}}$ are the differentials of the associated double complexes, the relative transgression map satisfies the following identity:
\begin{equation}\label{RelativeTransgressionBoundary}
	\psi_{\mathfrak{U}, f} \circ D_{\mathfrak{U}} = (-1)^{m}D_{\mathfrak{W}} \circ \psi_{\mathfrak{U}, f}.
\end{equation}
\end{Lemma}
\paragraph{Proof:} From \eqref{Cone}, \eqref{TransgressionBoundary} and \eqref{DiagramM} we get:
	\[\begin{split}
	\psi_{\mathfrak{U}, f} \circ D_{\mathfrak{U}}(\alpha_{1}, f_{*}\alpha_{2}) &= \psi_{\mathfrak{U}, f}(D_{\mathfrak{U}}\alpha_{1}, \varphi_{f,p,p-1}(\alpha_{1}) - f_{*}D_{f^{-1}\mathfrak{U}}\alpha_{2}) \\
	&= \pi_{1}^{*}\psi_{\mathfrak{U}}(D_{\mathfrak{U}}\alpha_{1}) - \pi_{2}^{*}\psi_{f^{-1}\mathfrak{U}}(f^{*}\varphi_{f,p,p-1}(\alpha_{1}) - D_{f^{-1}\mathfrak{U}}(\alpha_{2})) \\
	&= \pi_{1}^{*}[(-1)^{m}D_{\mathfrak{V}}\psi_{\mathfrak{U}}(\alpha_{1}) + r^{*}(\partial \psi)_{\mathfrak{U}}(\alpha_{1})] \\
	&\phantom{XXXXXX} - \pi_{2}^{*}[l_{f}^{*}(\partial \psi)_{\mathfrak{U}}(\alpha_{1}) - (-1)^{m-1}D_{\mathfrak{V}'}\psi_{f^{-1}\mathfrak{U}}(\alpha_{2})] \\
	&= (-1)^{m}\pi_{1}^{*}D_{\mathfrak{V}}\psi_{\mathfrak{U}}(\alpha_{1}) + (-1)^{m-1}\pi_{2}^{*}D_{\mathfrak{V}'}\psi_{f^{-1}\mathfrak{U}}(\alpha_{2}) \\
	&= (-1)^{m}D_{\mathfrak{W}} \circ \psi_{\mathfrak{U}, f}(\alpha_{1}, f_{*}\alpha_{2}).
\end{split}\]
$\square$ \\

It follows that the map \eqref{RelativeTransgression} defines a map in cohomology:
\begin{equation}\label{RelativeTransgressionCohomologyCover}
	\psi_{\mathfrak{U}, f}: H^{k}(\mathfrak{U}, S_{X,p}, S_{Y,p-1}, f) \rightarrow H^{k-m}(\mathfrak{W}, S_{\mathcal{M}_{f},p-m}).
\end{equation}
Using the same argument of \cite[Proposition 2.3]{GT2}, one can prove that:
\begin{Theorem} Considering the direct limit with respect to the covers of $X$, from \eqref{RelativeTransgressionCohomologyCover} we get a well-defined \emph{transgression map in relative Deligne cohomology}:
\begin{equation}\label{RelativeTransgressionCohomology}
	\psi_{f}: H^{k}(S_{X,p}, S_{Y,p-1}, f) \rightarrow H^{k-m}(S_{\mathcal{M}_{f},p-m}).
\end{equation}
$\square$
\end{Theorem}
When $m = k = p$, $\psi_{f}(\alpha_{1}, f_{*}\alpha_{2})$ is a function $\psi_{f}(\alpha_{1}, f_{*}\alpha_{2}): \mathcal{M}_{f} \rightarrow U(1)$, which, evaluated on $(g, g') \in \mathcal{M}$, corresponds to the holonomy of $(\alpha_{1}, f_{*}\alpha_{2})$ over the relative cycle $(g_{*}[M], -g'_{*}[\partial M])$, with $[M]$ and $[\partial M]$ computed via any triangulation of $M$. We now generalize the definition of holonomy to any cycle, not necessarily of the form $(g_{*}[M], -g'_{*}[\partial M])$ for a manifold $M$. \\

Let us consider formula \eqref{Transgression} for $M = \Delta^{m} = \{x \in \mathbb{R}^{m+1}: 0 \leq x_{i} \leq 1, x_{1} + \cdots + x_{m+1} = 1\}$. We fix a singular $m$-chain $\sigma = \sum_{i \leq a} n_{i} \sigma_{i}$ for $\sigma_{i}: \Delta^{m} \rightarrow X$. We remark that, for the moment, we have put an upper-bound $a$ on the number of simplices with possibly non-zero coefficient: we call $C_{m}^{[a]}(X)$ the set of such simplices. Given a cover $\mathfrak{U}$ on $X$, let $\mathfrak{V}$ be the cover induced on $C^{\infty}(\Delta^{m}, X)$ by $\mathfrak{U}$.\footnote{We cannot consider only the trivial triangulation of $\Delta^{m}$ because its image under a simplex is not necessarily contained in an element of $\mathfrak{U}$. That's why also the cover $\mathfrak{V}$ can be non-trivial.} We consider the projections $\pi_{i}: C_{m}^{[a]}(X) \rightarrow C^{\infty}(\Delta^{m}, X)$ for $i = 1, \ldots, a$, and we choose a common refinement $\mathfrak{V}'$ of the pull-backs $\pi_{i}^{-1}\mathfrak{V}$. We can now apply formula \eqref{Transgression} to each simplex and obtain a map:
\begin{equation}\label{PsiChains}
\begin{split}
	\psi_{\mathfrak{U}}^{[a]}: \;&C^{k}(\mathfrak{U}, S_{X,p}) \rightarrow C^{k-m}\bigl(\mathfrak{V}', S_{C_{m}^{[a]}(X),p-m}\bigr) \\
	& \psi_{\mathfrak{U}}^{[a]}(\alpha) = \sum_{i=1}^{a} n_{i} \cdot \pi_{i}^{*}\psi_{\mathfrak{U}}(\alpha).
\end{split}
\end{equation}
We also consider the natural projections $\pi'_{i} = r \circ \pi_{i}: C_{m}^{[a]}(X) \rightarrow C^{\infty}(\partial \Delta^{m}, X)$, and the boundary $\partial: C_{m}^{[a]}(X) \rightarrow C_{m-1}^{[a(m+1)]}(X)$. By linearity the map \eqref{PsiChains} still satisfies \eqref{TransgressionBoundary}, replacing $r^{*}$ with $(-\partial)^{*}$:
\begin{equation}\label{TransgressionBoundaryChains}
\begin{split}
	\psi_{\mathfrak{U}}^{[a]} \circ D_{\mathfrak{U}}(\alpha) &= \sum_{i=1}^{a} n_{i} \cdot \pi_{i}^{*}\psi_{\mathfrak{U}}D_{\mathfrak{U}}(\alpha) = \sum_{i=1}^{a} n_{i} \bigl[ (-1)^{m} \pi_{i}^{*}(D_{\mathfrak{V}} \psi_{\mathfrak{U}}(\alpha) ) + \pi_{i}^{*}r^{*} (\partial \psi)_{\mathfrak{U}}(\alpha) \bigr] \\
	&= (-1)^{m}\Bigl(D_{\mathfrak{V}} \sum_{i=1}^{a} n_{i} \cdot \pi_{i}^{*}\psi_{\mathfrak{U}}(\alpha) \Bigr) + \sum_{i=1}^{a} n_{i} \cdot {\pi_{i}'}^{*}(\partial\psi)_{\mathfrak{U}}(\alpha) \\
	&= (-1)^{m} D_{\mathfrak{V}}\psi_{\mathfrak{U}}^{[a]}(\alpha) - \partial^{*}\psi_{\mathfrak{U}}^{[a(m+1)]}(\alpha).
\end{split}
\end{equation}
Given a map $f: Y \rightarrow X$, we get a commutative diagram analogous to \eqref{DiagramM} in the following way:
\begin{equation}\label{DiagramMChains}
	\xymatrix{
	& C_{p}(X) \ar[dr]^{-\partial} & \\
	Z_{p}(f) \ar[ur]^{\pi_{1}} \ar[dr]_{\pi_{2}} & & C_{p-1}(X). \\
	& C_{p-1}(Y) \ar[ur]_{f_{*}}
}
\end{equation}
Therefore, we can define in the same way the map \eqref{RelativeTransgression} and, with the same proof of lemma \ref{RelativeTransgressionBoundaryThm}, for $k = m = p$ we get a map analogous to \eqref{RelativeTransgressionCohomology}:
	\[\Hol_{f}^{[a]}: H^{p}(S_{X,p}, S_{Y,p-1}, f) \rightarrow \Gamma(\underline{U}(1)_{Z_{p}^{[a]}(f)}).
\]
Since this function do not depend on the cover any more, we can take the direct limit on $a$ and define:
\begin{equation}\label{HolonomyMap}
	\Hol_{f}: H^{p}(S_{X,p}, S_{Y,p-1}, f) \rightarrow \Gamma(\underline{U}(1)_{Z_{p}(f)}).
\end{equation}
The function \eqref{HolonomyMap} is the relative holonomy. We now prove that the map \eqref{RelativeTransgressionCohomology} defines a natural isomorphism between the classes in $H^{p}(S_{X,p}, S_{Y,p-1}, f)$ and the relative Cheeger-Simons character of type II of degree $p$.
\begin{Theorem} There is a natural isomorphism:
\begin{equation}\label{Xip}
	\Xi^{p}: H^{p}(S_{X,p}, S_{Y,p-1}, f) \overset{\simeq}\longrightarrow CS_{II}^{p}(f)
\end{equation}
defined by $\Xi^{p}[(\alpha_{1}, f_{*}\alpha_{2})] := (\chi, \omega, \rho)$, where $\omega$ is the curvature of $[\alpha_{1}]$, $\rho$ the form induced by the strong topological trivialization on $Y$ and, for $(C, C') \in Z_{p}(f)$:
	\[\chi(C, C') = \frac{1}{2\pi i} \log \Hol_{f}[(\alpha_{1}, f_{*}\alpha_{2})](C, C') \mod \mathbb{Z},
\]
where $\Hol_{f}$ is defined by \eqref{HolonomyMap}.
\end{Theorem}
\paragraph{Proof:} We start showing that $\Xi^{p}$ is well-defined. If $(C, C') = \partial(D, D')$, then, we think of $(\alpha_{1}, f_{*}\alpha_{2})$ as a cochain in $C^{p+1}(S_{X,p}, S_{Y,p-1})$, so that $D_{\mathfrak{U}}(\alpha_{1}, f_{*}\alpha_{2}) = ((1, 0, \ldots, 0, \omega), f_{*}(1, 0, \ldots, 0, \rho))$, for $\omega$ the curvature of $\alpha_{1}$ and $\rho$ the connection determined by the topological trivialization $\alpha_{2}$. From \eqref{Cone} and \eqref{TransgressionBoundaryChains} we get:
	\[\begin{split}
	\psi_{\mathfrak{U}, f} \circ D_{\mathfrak{U}}(&\alpha_{1}, f_{*}\alpha_{2}) = \psi_{\mathfrak{U}, f}(D_{\mathfrak{U}}\alpha_{1}, \varphi_{f,p,p-1}(\alpha_{1}) - f_{*}D_{f^{-1}\mathfrak{U}}\alpha_{2}) \\
	&= \pi_{1}^{*}\psi_{\mathfrak{U}}(D_{\mathfrak{U}}\alpha_{1}) - \pi_{2}^{*}\psi_{f^{-1}\mathfrak{U}}(f^{*}\varphi_{f,p,p-1}(\alpha_{1}) - D_{f^{-1}\mathfrak{U}}(\alpha_{2})) \\
	&= \pi_{1}^{*}[(-1)^{m}D_{\mathfrak{V}}\psi_{\mathfrak{U}}(\alpha_{1}) - \partial_{X}^{*}\psi_{\mathfrak{U}}(\alpha_{1})] \\
	&\phantom{XXXX} - \pi_{2}^{*}[(f_{*})^{*}\psi_{\mathfrak{U}}(\alpha_{1}) - (-1)^{m-1}D_{\mathfrak{V}'}\psi_{f^{-1}\mathfrak{U}}(\alpha_{2}) + \partial_{Y}^{*}\psi_{f^{-1}\mathfrak{U}}(\alpha_{2})] \\
	&= (-1)^{m}D_{\mathfrak{W}} \circ \psi_{\mathfrak{U}, f}(\alpha_{1}, f_{*}\alpha_{2}) - \pi_{1}^{*}\partial_{X}^{*} \psi_{\mathfrak{U}}(\alpha_{1}) \\
	&\phantom{XXXX} - \pi_{2}^{*}(f_{*})^{*}\psi_{\mathfrak{U}}(\alpha_{1}) - \pi_{2}^{*}\partial_{Y}^{*} \psi_{f^{-1}\mathfrak{U}}(\alpha_{2}).
\end{split}\]
From the l.h.s.\ we get:
	\[\begin{split}
	\psi_{\mathfrak{U}, f}((1, &0, \ldots, 0, \omega), (1, 0, \ldots, 0, \rho))(D, -D') = \\ 
	& \psi_{\mathfrak{U}}(1, 0, \ldots, 0, \omega)(D) - \psi_{f^{-1}\mathfrak{U}} (1, 0, \ldots, 0, \rho)(D') = \int_{D} \omega + \int_{D'} \rho,
\end{split}\]
the last equality following directly from the definition of $\psi_{\mathfrak{U}}$ \cite[Definition 2.1]{GT2}. In the r.h.s., $\psi_{\mathfrak{U}, f}(\alpha_{1}, f_{*}\alpha_{2})$ is a cochain of degree $-1$, therefore vanishing. What remains evaluating in $(D, -D')$ is:
	\[\begin{split}
	-\partial_{X}^{*} \psi_{\mathfrak{U}}&(\alpha_{1})(D) - (f_{*})^{*}\psi_{\mathfrak{U}}(\alpha_{1})(-D') - \partial_{Y}^{*} \psi_{f^{-1}\mathfrak{U}}(\alpha_{2})(-D') \\
	& = \log[\Hol(\alpha_{1})(\partial D) - \Hol(f^{*}\alpha_{1})(-D') + \Hol(\alpha_{2}(-\partial D')]\\
	& = \log\Hol_{f}(\alpha_{1}, f_{*}\alpha_{2})(\partial D + f_{*}D', -\partial D') = \log\Hol_{f}(\alpha_{1}, f_{*}\alpha_{2})(C, C').
\end{split}\]
Therefore we get formula \eqref{RelativeCS2}. In order to see that \eqref{Xip} is an isomorphism, we can construct an exact sequence analogous to \eqref{LongDeligne2} in the following way:
\begin{equation}\label{LongDeligne2CS}
\begin{split}
	\cdots & \longrightarrow H^{p-1}(X, Y; \mathbb{R}/\mathbb{Z}) \longrightarrow H^{p-1}(X; \mathbb{R}/\mathbb{Z}) \longrightarrow CS^{p-1}(Y) \\
	& \longrightarrow CS_{II}^{p}(X,Y) \longrightarrow CS^{p}(X) \longrightarrow H^{p+1}(Y; \mathbb{Z}) \\
	& \longrightarrow H^{p+2}(X, Y; \mathbb{Z}) \longrightarrow H^{p+2}(X; \mathbb{Z}) \longrightarrow H^{p+2}(Y; \mathbb{Z}) \longrightarrow \cdots.
\end{split}
\end{equation}
There is a natural map from \eqref{LongDeligne2} to \eqref{LongDeligne2CS} defined via \eqref{Xip} and the analogous morphisms $H^{p-1}(S_{Y,p-1}) \rightarrow CS^{p-1}(Y)$ and $H^{p}(S_{X,p}) \rightarrow CS^{p}(X)$; since the latter are isomorphisms, it follows from the five lemma that even \eqref{Xip} is. $\square$ \\

We can now extend the statement to relative Cheeger-Simons character of any type.
\begin{Theorem} We call $CS_{t}^{p}(X,Y)$ the group of relative Cheeger-Simons characters of type $t$ of degree $p$. There are natural isomorphisms for $t = 0,1,2$:
\begin{equation}\label{Xipt}
	\Xi^{p}_{t}: H^{p}(S_{X,p}, S_{Y,p-t}) \overset{\simeq}\longrightarrow CS_{t+1}^{p}(X,Y).
\end{equation}
\end{Theorem}
\paragraph{Proof:} For $t = 1$ the isomorphism is \eqref{Xip}. For $t = 0$, it is enough to consider the inclusions \eqref{MapPP1} and $(\chi, \omega) \rightarrow (\chi, \omega, 0)$ and apply again \eqref{Xip}. For type III, given $\alpha \in H^{p}(S_{X,p}, S_{Y,0})$, we choose a counterimage $\tilde{\alpha} \in H^{p}(S_{X,p}, S_{Y,p-1})$ via the surjective map \eqref{MapP10}, we apply $\eqref{Xip}$ to it, and we apply to $\Xi^{p}(\tilde{\alpha})$ the surjective map from type II characters to type III ones defined via \eqref{RelativeCS3}. $\square$

\begin{Theorem}\label{ThmXip3} There is a natural isomorphism:
\begin{equation}\label{Xip3}
	\Xi^{p}_{3}: \overline{H}^{p}(S_{X,p}, S_{Y,p}) \overset{\simeq}\longrightarrow CS_{IV}^{p}(X,Y).
\end{equation}
\end{Theorem}
\paragraph{Proof:} The numerator of \eqref{HBar} corresponds under the isomorphism \eqref{Xip} to $CS_{II'}(f)$, because the latter is exactly the group of characters vanishing on $Y$. For the denominator, because of \eqref{LongDeligne2} the kernel of the map $H^{p-1}(S_{Y,p-1}) \rightarrow H^{p}(S_{X,p}, S_{Y,p-1}, f)$ is made by the image of flat classes on $X$, therefore, when composing with $H^{p-1}(S_{X,p-1})$ $\rightarrow H^{p-1}(S_{Y,p-1})$, only the curvature of the original class in $X$ is meaningful. In particular, given $\alpha_{1} \in H^{p-1}(S_{X,p-1})$, the maps $H^{p-1}(S_{X,p-1}) \rightarrow H^{p-1}(S_{Y,p-1}) \rightarrow H^{p}(S_{X,p}, S_{Y,p-1}, f)$ act as $\alpha_{1} \rightarrow \varphi_{f,p,p-1}\alpha_{1} \rightarrow (0, \varphi_{f,p,p-1}\alpha_{1})$, therefore, when we compute the holonomy on $(C,C')$, we get $-\Hol(\varphi_{f,p,p-1}\alpha_{1})(C') = -\Hol(\alpha_{1})(f_{*}C') = \Hol(\alpha_{1})(\partial C) = \exp 2\pi i\int_{C} \omega$, for $\omega$ the curvature of $\alpha_{1}$. Hence we get exactly the image of \eqref{phif}. $\square$ \\

Finally, let us consider a fiber bundle $E \rightarrow X$, where $X$ and the fiber $M$ are smooth manifolds (even with boundary), $M$ is oriented and $m = \dim\,M$. Given covers $\mathfrak{U} = \{U_{i}\}_{i \in I}$ of $X$ and $\mathfrak{U} = \{V_{i}\}_{i \in J}$ of $M$, a triangulation $K$ of $M$ and a map $\phi: K \rightarrow J$, in \cite[Definition 2.1]{GT} the authors define a cover $\mathfrak{W}$ of $E$ naturally associated to $\mathfrak{U}$, and a transgression map
\begin{equation}\label{IntegrationMap}
	\psi_{K, \phi, \mathfrak{U}}: C^{k}(\mathfrak{W}, S_{E,p}) \rightarrow C^{k-m}(\mathfrak{U}, S_{X,p-m}).
\end{equation}
Let $E'$ be the sub-bundle of $E$ obtained considering the boundary of each fiber, $r: E \rightarrow E'$ be the restriction map, and $(\partial \psi)_{\partial K, \partial \phi, \mathfrak{U}}$ the map \eqref{IntegrationMap} obtained considering $E'$ instead of $E$. If $D_{\mathfrak{U}}$ and $D_{\mathfrak{W}}$ are the differentials of the associated double complexes, then \eqref{IntegrationMap} satisfies \cite[Theorem 3.1]{GT}
\begin{equation}\label{IntegrationMapBoundary}
	\psi_{K, \phi, \mathfrak{U}} \circ D_{\mathfrak{W}} = (-1)^{m}D_{\mathfrak{U}} \circ \psi_{K, \phi, \mathfrak{U}} + r^{*} \circ (\partial \psi)_{\partial K, \partial \phi, \mathfrak{U}}.
\end{equation}
\begin{Def} With the data defined above, for $\iota: E' \hookrightarrow E$ the embedding, we define the \emph{relative integration map}
\begin{equation}\label{RelativeIntegration}
	\psi'_{K, \phi, \mathfrak{U}}: C^{k}(\mathfrak{W}, S_{E,p}, S_{E',p-1}) \rightarrow C^{k-m}(\mathfrak{U}, S_{X,p-m}).
\end{equation}
as:
	\[\psi'_{\mathfrak{U}}(\alpha_{1}, \iota_{*}\alpha_{2}) := \psi_{\mathfrak{U}}(\alpha_{1}) - r^{*}(\partial \psi)_{\partial K, \partial \phi, \mathfrak{U}\vert_{Y}}(\alpha_{2}).
\]
\end{Def}
\begin{Theorem}\label{RelativeIntegrationBoundaryThm} If $D_{\mathfrak{U}}$ and $D_{\mathfrak{V}}$ are the differentials of the associated double complexes, the relative transgression map satisfies the following identity:
\begin{equation}\label{RelativeIntegrationBoundary}
	\psi'_{K, \phi, \mathfrak{U}} \circ D_{\mathfrak{W}} = (-1)^{m}D_{\mathfrak{U}} \circ \psi'_{K, \phi, \mathfrak{U}}.
\end{equation}
\end{Theorem}
\paragraph{Proof:} We argue as in lemma \ref{RelativeTransgressionBoundaryThm}. $\square$ \\

It follows that the map \eqref{RelativeIntegration} defines a map in cohomology $\psi'_{K, \phi, \mathfrak{U}}: H^{k}(\mathfrak{W}, S_{E,p}, S_{E',p-1}) \rightarrow H^{k-m}(\mathfrak{U}, S_{X,p-m})$. Using the same argument of \cite[Proposition 2.3]{GT2} one can prove that:
\begin{Theorem} Considering the direct limit with respect to the covers of $X$, from \eqref{RelativeIntegration} we get a well-defined map:
\begin{equation}\label{RelativeIntegrationCohomology}
	\psi': H^{k}(S_{E,p}, S_{E',p-1}) \rightarrow H^{k-m}(S_{X,p-m}).
\end{equation}
$\square$
\end{Theorem}

\section*{Acknowledgements}

The author is financially supported by FAPESP (Funda\c{c}\~ao de Amparo \`a Pesquisa do Estado de S\~ao Paulo).


\end{document}